# ROC AND THE BOUNDS ON TAIL PROBABILITIES VIA THEOREMS OF DUBINS AND F. RIESZ

By Eric Clarkson,[1] J. L. Denny and Larry Shepp[2]

*University of Arizona, University of Arizona and Rutgers University*

For independent $X$ and $Y$ in the inequality $P(X \leq Y + \mu)$, we give sharp lower bounds for unimodal distributions having finite variance, and sharp upper bounds assuming symmetric densities bounded by a finite constant. The lower bounds depend on a result of Dubins about extreme points and the upper bounds depend on a symmetric rearrangement theorem of F. Riesz. The inequality was motivated by medical imaging: find bounds on the area under the Receiver Operating Characteristic curve (ROC).

**1. Introduction.** We give sharp upper and lower bounds on

$$P(X \leq Y + \mu),$$

where the independent variables $X$ and $Y$ have zero means and satisfy either unimodality or symmetry conditions. The lower bounds assume unimodality and use a theorem of Dubins [5] about extreme points, while the upper bounds assume symmetry and use a theorem of F. Riesz [13] about symmetric rearrangements. We emphasize that our basic inequalities in the lower-bound case are known, proved earlier by various authors starting with Gauss. Our justification for proving the known theorems is mainly to show that the bounds are sharp and perhaps to indicate another approach.

Both bounds were motivated by a widely used methodology in medical imaging, the ROC (Receiver Operating Characteristic) curve, known to statisticians as a power function (not necessarily of the most powerful test). A widely used interpretation of the ROC curve is the AUC (Area under the Curve), defined next.

Received September 2007; revised April 2008.
[1]Supported by NIH Grants P41 EB 002305 and R37 EB00803.
[2]Supported in part by NSF Grant 0203086—Statistical Methods in Fast Functional MRI.

Thus, we are given variables $X_i$, $i = 0, 1$ having continuous $F_i(x) = P(X_i \leq x)$. Letting for $0 < \alpha < 1$, $x(\alpha) = \sup\{x : 1 - F_0(x) \geq \alpha\}$, the ROC curve is

$$\alpha \longmapsto x(\alpha),$$

the AUC being

$$AUC(X_0, X_1) = \int_0^1 (1 - F_1(x(\alpha))) \, d\alpha. \tag{1}$$

If the variables $X_i$ are independent, then (1) equals

$$P(X_0 \leq X_1), \tag{2}$$

an identity first proved by Bamber [1], at least in the medical imaging literature.

This is a role of the AUC. An experimenter wants to compare two medical imaging modalities to decide which best detects a tumor. For example, one may compare X-rays against MRI images, although often one compares "filters" for the same modality [other modalities include positron emission tomography (PET), single-photon computed emission tomography (SPECT) and ultrasound]. Imagine that a large number of experimenters each test the hypothesis of $F_0$ (no tumor) against $F_1$ (tumor) and that each chooses a level of significance $\alpha$ according to a uniform distribution. (The use of a random $\alpha$ reflects the differing levels of significance of different experimenters.) For each experimenter, the hypothesis $F_0$ is rejected in favor of $F_1$ if a scalar observable exceeds a constant $x(\alpha)$. Thus the AUC gives the average of the power function $1 - F_1(x(\alpha))$.

The equality of (2) and (1) for continuous distributions leads to a second widely used method, the 2AFC (Two-Alternative Forced Choice). In this case the experimenter is confronted with two choices, perhaps two different imaging modalities, perhaps a "signal" or "no signal." The experimenter uses a test statistic, rejecting the hypothesis of no signal in favor of a signal if the test statistic is large. More precisely, the experimenter computes the test statistic, applies it to the two data sets, the signal and the nonsignal (not knowing which is which), and chooses the data set giving the larger value of the statistic as the signal. In (2), the distribution of $X_1$ is that of the statistic when the signal is present, $X_0$ when the signal is absent.

The ROC was developed during World War II for analyzing the performance of radar systems. Today ROC analysis is regularly used in the health care industry and by the Federal Drug Administration to evaluate new imaging systems, diagnostic tests, treatments and pharmaceuticals. Often, if not invariably, a Gaussian assumption is made on a test statistic, typically the log-likelihood ratio. The bounds obtained here are of course weaker. For example, without assuming a Gaussian distribution or the equality of modes (but assuming $\mu_X = \mu_Y$, $\sigma_X = \sigma_Y = 1$), the unimodality lower bound when $\mu = 2\sqrt{6} \cong 4.899$ is

$$P(X \leq Y + 2\sqrt{6}) \geq 0.5,$$

while a Gaussian distribution gives a lower bound greater than 0.99966 (since $\sqrt{12} \cong 3.4641$); claimed differences may be the result of the Gaussian assumption. Some sort of compromise is needed.

ROC analysis in medical imaging has an enormous literature. We mention the book of Swets and Pickett [14] and the papers of Metz [11], Clarkson [3] and Barrett, Abbey and Clarkson [2]. We also mention that researchers in medicine and in psychophysics—the branch of psychology that studies the relations between physical stimuli and sensory response—use a functional relation between the AUC and the SNR (Signal-to-Noise Ratio), where in medical imaging, the SNR is defined as the ratio of the mean pixel value to their standard deviation.

In this paper we study the behavior of (2) under unimodality or symmetry assumptions. Lower and upper confidence bounds on a translation parameter $\mu$ defined below are clearly available although not discussed here. Because (2) and (1) are equal for continuous distributions, we state our results in terms of (2). For the lower bounds in the unimodal case we constrain the variances to equal 1. For the upper bounds in the symmetric case we constrain the densities to be bounded by $b < \infty$. Although the second constraint may appear unfamiliar (if not unnatural), it is easy to see that neither constraint is relevant to the other case.

**2. Dubins's theorem—lower bounds for $P(X \leq Y + \mu)$.** Throughout this section we assume that $X_0$ has a mode equal to zero. Rather than assuming that $X_1$ has a mode equal to $\mu > 0$, we find it convenient to assume that $X_1$ has a mode equal to zero and then study $X_1 + \mu$.

We begin with sharp *upper* bounds for *symmetric unimodal* distributions and we recall the definition of unimodal distributions (page 155 of Feller [6]): a distribution function $F$ is *unimodal at* $m$ if $F$ is convex on $(-\infty, m)$ and concave on $(m, \infty)$. Note that $F$ may assign positive mass to the point $m$ and that a mode need not be unique.

Gauss (see Pukelsheim [12]) proved the following inequality for a variable $X$ having a continuous unimodal distribution, where $\tau^2 = E[(X-m)^2]$:

$$P(|X - m| > s) \leq \begin{cases} 1 - (s/\sqrt{3}\tau), & s \leq 2\sqrt{3}\tau, \\ 4\tau^2/9s^2, & s \geq 2\sqrt{3}\tau. \end{cases}$$

Pukelsheim [12] gives a useful survey of the inequalities which followed that of Gauss. The method of extreme points was used earlier by Dharmadhikari and Joag-Dev [4].

We need a special case of a theorem of L. Dubins [5].

THEOREM 1 (Dubins). *Let $A$ be a compact convex subset of a locally convex space and let $T$ be a real continuous linear functional on $A$. Then each extreme point of $A \cap \{x : T(x) = y\}$ is a convex combination of at most two extreme points of $A$.*

The inequality in the next lemma and its corollary are known, mentioned in the references above. What is new (we believe) is a sharp one-sided inequality. Dubins's theorem leads naturally to the bounds and the distributions achieving the bounds.

In order to use a compactness argument on the space of distributions, we initially assume that all distributions are supported on the interval $[-N, N]$. Since $N$ is arbitrary, it is easy to verify that all assertions will extend to the case that the distributions are supported on $(-\infty, \infty)$.

LEMMA 2. *Fix $t > 0$ and assume that $t \leq 2N/3$. If $X$ has a symmetric unimodal distribution supported on $[-N, N]$ and $\operatorname{var}(X) = 1$, then*

$$P(X \geq t) \leq \begin{cases} 1/2 - t/(2\sqrt{3}), & 0 < t \leq 2/\sqrt{3}, \\ 2/(9t^2), & 2/\sqrt{3} < t \leq 2N/3. \end{cases}$$

*If $0 < t \leq 2/\sqrt{3}$, the bound is obtained at the density*

$$f(x) = \begin{cases} 1/(2\sqrt{3}), & |x| < \sqrt{3}, \\ 0, & |x| \geq \sqrt{3}. \end{cases}$$

*If $2/\sqrt{3} < t \leq 2N/3$, the bound is obtained at the distribution*

$$(1 - \beta_t)\delta_0 + \beta_t f_t(x),$$

*where $\beta_t = 4/(3t^2)$ and*

$$f_t(x) = \begin{cases} 1/(3t), & |x| < 3t/2, \\ 0, & |x| \geq 3t/2. \end{cases}$$

PROOF. We fix $N \geq 2$, let $\mathcal{F}$ denote the set of *symmetric unimodal* distributions on $[-N, N]$, and note that $\mathcal{F}$ is a compact convex set in the weak topology. It is easy to see that the extreme points of $\mathcal{F}$ are the Dirac probability $\delta_0$ and the boxcar densities $f_a(x)$, where for $0 < a \leq N$,

$$f_a(x) = \begin{cases} 1/(2a), & |x| < a, \\ 0, & |x| \geq a. \end{cases}$$

Define a continuous linear functional $T$ on $\mathcal{F}$ by

$$T(F) = \int_{-N}^{N} x^2 \, dF(x)$$

and let $\mathcal{F}(1) \subset \mathcal{F}$ denote the compact convex set where $T(F) = 1$, that is, the variance equals 1. Dubins's theorem implies that the extreme points of $\mathcal{F}(1)$ are the distributions of the form

(3) $$(1-\beta)\delta_0 + \beta F_2,$$

(4) $$(1-\beta)F_1 + \beta F_2,$$

where the $F_i$ are distribution functions having boxcar densities and $0 < \beta < 1$, and the variance of (3) and of (4) equals 1. Fix $0 < t \leq 2N/3$ and define the continuous linear functional on $\mathcal{F}(1)$,

$$T_t(F) \to \int_t^N dF(x).$$

Because the maximum of $T_t(F)$ on $\mathcal{F}(1)$ is obtained at an extreme point, to prove the theorem it suffices to calculate the values of $T_t(F)$ on (3) and (4). □

Now let $X$ have a unimodal distribution function $H$ with mode at the origin, $var(X) = 1$, mean $\mu_X$, and let $X^s$ denote the variable with the *symmetric unimodal* distribution function $1/2(H(x) + 1 - H(-x))$. Then $X^s$ has variance $1 + \mu_X^2$.

COROLLARY 3. *For $t > 0$,*

$$P(|X| > t(1+\mu_X^2)^{1/2}) = 2P(X^s > t(1+\mu_X^2)^{1/2})$$
$$\leq \begin{cases} 1 - t/\sqrt{3}, & t < 2/\sqrt{3}, \\ 4/(9t^2), & t \geq 2/\sqrt{3}. \end{cases}$$

THEOREM 4. *Let $X$ have a unimodal distribution with a mode at the origin, $var(X) = 1$, and mean $\mu_X$. Then for $t > 0$,*

$$P(X > t) \leq \begin{cases} 1 - t/(2\sqrt{3}), & 0 < t \leq 4/\sqrt{3} \cong 2.3094, \\ 4(1+\mu_X^2), & (9t^2)4/\sqrt{3} < t. \end{cases}$$

*If $t \leq 4/\sqrt{3}$, then the bound is obtained at the density*

(5) $$f(x) = 1/(2\sqrt{3}), \qquad 0 < x < 2\sqrt{3}.$$

*Fix $t > 4/\sqrt{3}$. The bound is obtained at the distribution*

$$\left(1 - \frac{4}{3u^2}\right)\delta_0(x) + \frac{8(u^2-1)^{1/2}}{9u^4} I_{A(u)}(x),$$

*where*

$$A(u) = \left(0, \frac{3u^2}{2(u^2-1)^{1/2}}\right),$$

$$\mu_X = \frac{1}{(u^2-1)^{1/2}},$$

*and $u$ satisfies: $t$ is the positive root of*

$$\text{(6)} \qquad t^3 - \frac{3u^2}{2(u^2-1)^{1/2}}t^2 + \frac{1}{2}\left(\frac{u^4}{u^2-1}\right)^{3/2} = 0.$$

REMARK 5. Always $u \geq 2/\sqrt{3}$ and when $u = 2/\sqrt{3}$ the positive root $t = 4/\sqrt{3}$. If $t$ is large so that $u$ is also large, then (6) is approximately the polynomial in the variable $t$,

$$t^3 - \frac{3ut^2}{2} + \frac{1}{2}u^3 = 0,$$

whose unique root $t$ satisfies $t = u$. Let $t$ be given and let $u(t)$ be such that $t$ is the root of (6) with $u(t)$ the constant. For example, if $t \cong 3.18198$, then $u(t) = 3.0$, and if $t \cong 8.063242$, then $u(t) = 8.0$. It is not difficult to verify that $t > u(t)$ and $\lim t/t(u) = 1$ as $u \to \infty$.

PROOF OF THEOREM 4. To find the upper bound, we claim that we may assume that unimodal $X \geq 0$. For if $X$ has a unimodal distribution with mode at the origin, $var(X) = 1$, we may define a unimodal $S \geq 0$ with mode at the origin, $var(S) = 1$, so that for all $t > 0$,

$$P(S > t) \geq P(X > t).$$

If $X \leq 0$ choose, say, $S = -X$. Otherwise, take that part of the distribution of $X$ which is supported on $(-\infty, 0)$ and put it on $\{0\}$; call the new variable $Y$. Because

$$E(Y^2) \leq E(X^2),$$
$$E(Y) \geq E(X),$$

we have

$$0 < \gamma^2 \equiv Var(Y) \leq 1.$$

Then $S \equiv Y/\gamma \geq 0$ is unimodal with mode at the origin, $var(S) = 1$, and for $t > 0$,

$$P(S > t) = P(X > \gamma t) \geq P(X > t).$$

Continuing the proof, from Corollary 3 and Lemma 2,

$$\text{(7)} \qquad P(X > t) = 2P(X^s > t) = 2P(X^s/(1+\mu_X^2)^{1/2} > t/(1+\mu_X^2)^{1/2})$$
$$\leq \begin{cases} 1 - t/(\sqrt{3}(1+\mu_X^2)^{1/2}), & 0 < t < 2(1+\mu_X^2)^{1/2}/\sqrt{3}, \\ 4(1+\mu_X^2)/(9t^2), & 2(1+\mu_X^2)^{1/2}/\sqrt{3} < t. \end{cases}$$

If $X$ has the density (5), then $\mu_X = \sqrt{3}$ and $X$ satisfies (7) for $0 < t \leq 4/\sqrt{3}$. For $u \geq 2/\sqrt{3}$ let $X$ have the distribution

$$\left(1 - \frac{4}{3u^2}\right)\delta_0(x) + \frac{8(u^2-1)^{1/2}}{9u^4} I_{(0, 3u^2/(2(u^2-1)^{1/2})}.$$

Then $\sigma_X = 1$ and since $\mu_X = 1/(u^2-1)^{1/2}, 1 + \mu_X^2 = u^2/(u^2-1)$. For $t \geq (2/\sqrt{3})(u^2/(u^2-1))^{1/2}$,

(8) $$P(Y > t) = \frac{4}{3u^2}\left(1 - \frac{2t(u^2-1)^{1/2}}{3u^2}\right).$$

Now fix $t > 4/\sqrt{3}$. We want a value of $u$ so that (8) equals

$$= \frac{4u^2}{9t^2(u^2-1)} \quad \left(= \frac{4(1+\mu_X^2)}{9t^2}\right).$$

This is satisfied by (6). □

Following Ibragimov [8], a unimodal distribution function is *strong unimodal* if its composition with any unimodal distribution function is unimodal. Ibragimov proved that a distribution function $F$ is strong unimodal if and only if $F$ is continuous unimodal and its density $f$ satisfies

$$x \to \ln(f(x))$$

is a concave function on the interior of the support of $F$.

Let independent $X$ and $Y$ have unimodal distributions with modes (not necessarily unique) $m_X$ and $m_Y$, means $\mu_X = \mu_Y = 0$ and standard deviations $\sigma_X = \sigma_Y = 1$. We assume that at least one of the unimodal distributions is strong unimodal and recall from [9] that a mode of $(X-Y)/\sqrt{2}$ satisfies

(9) $$|m_{(X-Y)/\sqrt{2}}| \leq \sqrt{3}.$$

COROLLARY 6. *Fix $\mu > \sqrt{6}$. Then*

(10) $$P(X \leq Y + \mu)$$
$$\geq \begin{cases} (\mu - \sqrt{6})/(2\sqrt{6}), & \sqrt{6} \leq \mu \leq \sqrt{6} + 4\sqrt{2/3} \cong 5.7155, \\ 1 - 32(3(\mu - \sqrt{6}))^{-2}, & \sqrt{6} + 4\sqrt{2/3} < \mu. \end{cases}$$

REMARK 7. The only case of interest, $P(X \leq Y + \mu) \geq 0.5$, requires $\mu \geq 2\sqrt{6} \cong 4.899$.

PROOF OF COROLLARY 6. If we define

$$Z = \tfrac{1}{\sqrt{2}}(X - Y - m_{X-Y}),$$

then $Z$ is unimodal with a mode at the origin, mean equal to $(-m_{X-Y})/\sqrt{2}$ and variance equal to 1. The left-hand side of (10) is

$$\text{(11)} \qquad P(Z \leq \tfrac{1}{\sqrt{2}}(\mu - m_{X-Y})).$$

Thus among all such unimodal variables $Z$, (11) is minimized with $m_{X-Y}/\sqrt{2} = \sqrt{3}$, using (9). To complete the proof it suffices to find a least upper bound for

$$P(Z \geq \tfrac{1}{\sqrt{2}}(\mu - \sqrt{6})), \qquad \mu > \sqrt{6}$$

using Theorem 4. □

**3. F. Rieszs theorem—upper bounds for $P(X \leq Y + \mu)$.** Recall the definition of the *symmetric rearrangement* of the indicator function of a Borel set ([10]; see also [7]). If $A \subset R$ is a Borel set of finite Lebesgue measure $\lambda(A)$, then the symmetric rearrangement of the set $A$, denoted by $A^*$, is the symmetric open interval so that $\lambda(A^*) = \lambda(A)$. We let the functions

$$I_A, \qquad I_A^*$$

denote the indicator functions of $A$ and $A^*$, respectively. The following is a special case of Riesz's theorem [13]. As in the previous section, we study distributions whose support is $[-N, N]$, arbitrary $N$.

THEOREM 8 (F. Riesz). *Let $A, B$ and $C$ be Borel sets of finite measure. Then*

$$\int_{-\infty}^{\infty}\int_{-\infty}^{\infty} I_A(x) I_B(x-y) I_C(y)\, dx\, dy \leq \int_{-\infty}^{\infty}\int_{-\infty}^{\infty} I_A^*(x) I_B^*(x-y) I_C^*(y)\, dx\, dy.$$

Fix $b > 0$ and $N \geq 1/(2b)$, and let $\mathcal{F}(b)$ denote the class of all *symmetric* distributions on $[-N, N]$ whose distribution functions satisfy a Lipschitz condition with Lipschitz constant $b$. Clearly $\mathcal{F}(b)$ is convex and the Lipschitz condition ensures that $\mathcal{F}(b)$ is an equicontinuous class; because $\mathcal{F}(b)$ is closed (the sup norm topology), Ascoli's theorem implies that $\mathcal{F}(b)$ is compact. We let distributions $F, G \in \mathcal{F}(b)$ and let $H \in \mathcal{F}(b)$ be the distribution with density

$$\text{(12)} \qquad u(x) = \begin{cases} b, & |x| < 1/(2b), \\ 0, & |x| > 1/(2b). \end{cases}$$

With these distributions we associate independent variables

$$\text{(13)} \qquad X \sim F, \qquad Y \sim G; \qquad U_0, U_1 \sim H.$$

THEOREM 9. *Fix $\mu > 0$. If $b\mu < 1$, then*
$$P(X \leq Y + \mu) \leq P(U_0 \leq U_1 + \mu) = b\mu + \tfrac{1}{2}(1 - (b\mu)^2).$$
*The inequalities are strict unless $f = u$ a.e. If $b\mu \geq 1$, then*
$$P(X \leq Y + \mu) \leq P(U_0 \leq U_1 + \mu) = 1.$$

The proof of the theorem rests on the following lemmas. We thank the referee for observing that the argument of the next lemma extends to an arbitrary probability space without symmetry conditions.

LEMMA 10. *Assume that $N \geq 1/(2b)$. The extreme points of $\mathcal{F}(b)$ are the distributions with the densities (up to a set of measure zero)*

(14) $$bI_B(x) = \begin{cases} b, & x \in B, \\ 0, & x \notin B, \end{cases}$$

*where symmetric $B \subset [-N, N]$ and the Lebesgue measure $\lambda(B) = 1/b$.*

PROOF. Given $G$, let $g$ be the density of $G$. Suppose that $g$ does not satisfy (14) up to a set of measure zero: for an $\varepsilon > 0$ there is a symmetric set $A$
$$A = \{x : \varepsilon < g(x) < b - \varepsilon\}$$
so that $\lambda(A) > 0$. Choose disjoint symmetric subsets of $A_0, A_1 \subset A$, $\lambda(A_i) > 0$, and constants $\delta_0, \delta_1 > 0$ so that
$$\delta_0 \lambda(A_0) = \delta_1 \lambda(A_1),$$
and so that for $x \in A_0 \cup A_1, i = 0, 1 \pmod 2$,
$$g(x) - \delta_i \geq 0, \qquad g(x) + \delta_{i+1} \leq b.$$
Define for $i = 0, 1 \pmod 2$,
$$g_i(x) = \begin{cases} g(x), & x \in (A_0 \cup A_1)^c, \\ g(x) + (-1)^{i+1}\delta_0, & x \in A_0, \\ g(x) + (-1)^i \delta_1, & x \in A_1. \end{cases}$$
Then $G$ is not extremal: letting $G_i$ be the distribution functions of the $g_i$, we have distinct $G_i \in \mathcal{F}(b)$ and $G = \tfrac{1}{2}G_1 + \tfrac{1}{2}G_2$. □

Fix $\mu > 0$ and define the continuous bilinear functional on the convex compact $\mathcal{F}(b) \times \mathcal{F}(b)$,
$$T(F, G) = \int_{-N}^{N} dF(x) \left( \int_{(x-\mu) \vee (-N)}^{N} dG(y) \right)$$
$$= P(X \leq Y + \mu).$$

LEMMA 11. $\max\{T(F,G) : G, H \in \mathcal{F}(b)\}$ *is obtained at extreme points.*

PROOF. Recall that a continuous linear functional on a compact convex set obtains its maximum at an extreme point. Letting $(F_n, G_n)$ satisfy $T(F_n, F_n) \uparrow \sup\{T(F,G)\}$, one uses a compactness argument on subsequences to verify the lemma. □

The proof of the next lemma follows from Lemma 10 and the definition of symmetric rearrangement.

LEMMA 12. *If $F \in \mathcal{F}(b)$ is an extreme point with density $f$, then the symmetric rearrangement $f^* = u$ [see (12)].*

LEMMA 13. *Fix $\mu > 0$. Let $Y_0, Y_1$ have densities $f_0, f_1$ whose distributions are extreme points of $\mathcal{F}(b)$. Then*

$$P(|Y_0 - Y_1| < \mu) \leq P(|U_0 - U_1| < \mu).$$

PROOF. Note that $I_{(-\mu,\mu)} = I^*_{(-\mu,\mu)}$. Using the theorem of Riesz, the symmetry, and the preceding lemma,

$$\begin{aligned}
P(|Y_0 - Y_1| < \mu) &= \int_R \int_R I_{(-\mu,\mu)}(x) f_0(x-y) f_1(y) \, dy \, dx \\
&\leq \int_R \int_R I_{(-\mu,\mu)}(x) f_0^*(x-y) f_1^*(y) \, dy \, dx \\
&= \int_R \int_R I_{(-\mu,\mu)}(x) u(x-y) u(y) \, dy \, dx \\
&= P(|U_0 - U_1| < \mu). \quad \square
\end{aligned}$$

PROOF OF THEOREM 9. This follows from the lemmas and the symmetry assumption used in

$$P(Y_0 \leq Y_1 + \mu) = \tfrac{1}{2}(1 + P(|Y_1 - Y_0| < \mu)). \quad \square$$

**Acknowledgments.** Anirban Dasgupta and Harry Barrett contributed to this paper.

E. Clarkson
College of Optical Sciences
Department of Radiology
University of Arizona
Tucson, Arizona 85721
USA
E-mail: clarkson@radiology.arizona.edu

J. L. Denny
Department of Mathematics
Department of Radiology
University of Arizona
Tucson, Arizona 85721
USA
E-mail: denny@u.arizona.edu

L. Shepp
Statistics Department
Hill Center, Rutgers University
New Brunswick, New Jersey 08854
USA
E-mail: shepp@stat.rutgers.edu